\documentclass[a4,10pt,draft]{article}
\usepackage{amsmath}
\usepackage{amssymb}
\usepackage{amsthm}

\newtheorem*{kr}{Theorem A}
\theoremstyle{definition}
\newtheorem{remark}{Remark}
\newtheorem*{chu}{Remark}
\begin{document}
\title{Note on relations among multiple zeta-star values}
\author{Masahiro Igarashi}
\date{November 14, 2011}
\maketitle
\begin{abstract}
In this note, by using several ideas of other researchers, we derive several relations among multiple zeta-star values 
from the hypergeometric identities of C. Krattenthaler and T. Rivoal.
\end{abstract}
\section{Introduction}
The multiple zeta-star value (MZSV for short) is defined by the multiple series
\begin{equation*}
 \zeta^{\star}(k_1,\ldots,k_n) := \sum_{0< m_1{\le} \cdots {\le}m_n}\frac{1}{m_1^{k_1} \cdots m_n^{k_n}},
\end{equation*}
where $k_1,\ldots ,k_n\in\mathbb{Z}_{\ge1}$ and $k_n\ge2$. The case $n=2$ was studied by L. Euler in \cite{eu}, 
and the general case was introduced by M. E. Hoffman in \cite{ho}. 
In \cite{za}, D. Zagier pointed out the connection between multiple zeta values 
and several mathematical objects. 
Further, in \cite[p.~509]{za}, D. Zagier gave a conjecture for the dimension of the $\mathbb{Q}$-vector space 
generated by all multiple zeta values with given weight. 
MZSVs satisfy various relations 
(see, e.g., \cite{ako}, \cite{ao}, \cite{aow}, \cite{eu}, \cite{ho}, \cite{ikoo}, 
\cite{ko}, \cite{kaw}, \cite{mu}, \cite{ow}, \cite{oz}).
\par
In this note, by using several ideas of other researchers, we will derive several relations among MZSVs 
from the hypergeometric identities of C. Krattenthaler and T. Rivoal in \cite{kr}.
\par
First we recall the definition of the generalized hypergeometric series. 
The generalized hypergeometric series is defined by the power series
\begin{equation*}
{}_{p+1}F_p\left(
\begin{array}{c}
a_1,\ldots,a_{p+1}\\
b_1,\ldots,b_p
\end{array}
;z\right)
:=
\sum_{m=0}^{\infty}\frac{(a_1)_m\cdots(a_{p+1})_m}{(b_1)_m\cdots(b_p)_m}\frac{z^m}{m!},
\end{equation*}
where $p\in\mathbb{Z}_{\ge1}$, $z$, $a_i\in\mathbb{C}$ ($i=1,\ldots,p+1$), 
$b_j\in\mathbb{C}\setminus\mathbb{Z}_{\le0}$ ($j=1,\ldots,p$), and $(a)_m$ denotes 
the Pochhammer symbol defined by
\begin{equation*}
(a)_m = \left\{\begin{alignedat}{3}
                     &a(a+1)\cdots(a+m-1) &\quad&\text{if $m\in\mathbb{Z}_{\ge1}$},\\
                                   &1 &\quad&\text{if $m=0$}.
               \end{alignedat}
\right.
\end{equation*}
The above power series converges absolutely for all $z\in\mathbb{C}$ such that $|z|=1$ 
provided $\mathrm{Re}\,(b_1+\cdots+b_p-a_1-\cdots-a_{p+1})>0$.
\par
In \cite{kr}, C. Krattenthaler and T. Rivoal proved the following hypergeometric identities:
\begin{kr}[C. Krattenthaler and T. Rivoal {\cite[Proposition 1 (i) and (ii)]{kr}}]
\par
$(i)$ Let $s$ be a positive integer, and let $a$, $b_i$, $c_i$ $(i=1,\ldots,s+1)$ be complex numbers. 
Suppose that the complex numbers $a$, $b_i$, $c_i$ $(i=1,\ldots,s+1)$ satisfy the conditions
\begin{equation*}
\begin{aligned}
&1+a-b_i,\,1+a-c_i\notin\mathbb{Z}_{\le0}\,\,\,\, \textit{for} \,\,\, i=1,\ldots,s+1;\\
&\mathrm{Re}\left((2s+1)(a+1)-2\sum_{i=1}^{s+1}(b_i+c_i)\right)>0;\\
&\mathrm{Re}\left(\sum_{i=r}^{s+1}A_i(1+a-b_i-c_i)\right)>0
\end{aligned}
\end{equation*}
for $r=2,\ldots,s+1$, and all possible choices of $A_i=1$ or $2$ $(i=2,\ldots,s)$, $A_{s+1}=1$. 
$($For the details of the choices of $A_i$, see $[30]$.$)$ 
Then the following identity holds$:$
\begin{equation*}
\begin{aligned}
&{}_{2s+4}F_{2s+3}\left(
\begin{array}{c}
a, \frac{a}{2}+1, b_1, c_1, \ldots, b_{s+1}, c_{s+1}\\
\frac{a}{2}, 1+a-b_1, 1+a-c_1, \ldots, 1+a-b_{s+1}, 1+a-c_{s+1}
\end{array}
;-1\right)\\
=
&\frac{\Gamma(1+a-b_{s+1})\Gamma(1+a-c_{s+1})}{\Gamma(1+a)\Gamma(1+a-b_{s+1}-c_{s+1})}\\
&\times\sum_{l_1,\ldots,l_s\ge0}\prod_{i=1}^{s}
\frac{(1+a-b_i-c_i)_{l_i}(b_{i+1})_{l_1+\cdots+l_i}(c_{i+1})_{l_1+\cdots+l_i}}
{l_i!(1+a-b_i)_{l_1+\cdots+l_i}(1+a-c_i)_{l_1+\cdots+l_i}}.
\end{aligned}
\end{equation*}
\par
$(ii)$ Let $s$ be a positive integer, and let $a$, $b_i$ $(i=1,\ldots,s)$, $c_j$ $(j=0,1,\ldots,s)$ be complex numbers. 
Suppose that the complex numbers $a$, $b_i$ $(i=1,\ldots,s)$, $c_j$ $(j=0,1,\ldots,s)$ satisfy the conditions
\begin{equation*}
\begin{aligned}
&1+a-b_i,\,1+a-c_j\notin\mathbb{Z}_{\le0}\,\,\,\,\textit{for}\,\,\, i=1,\ldots,s,\,j=0,1,\ldots,s;\\
&\mathrm{Re}\left(2s(a+1)-2c_0-2\sum_{i=1}^{s}(b_i+c_i)\right)>0;\\
&\mathrm{Re}\left(\sum_{i=r}^{s}A_i(1+a-b_i-c_i)\right)>0,\\
&\mathrm{Re}\left(1+a-c_0-b_1-c_1+\sum_{i=2}^{s}A_i(1+a-b_i-c_i)\right)>0
\end{aligned}
\end{equation*}
for $r=2,\ldots,s$, and all possible choices of $A_i=1$ or $2$ $(i=2,\ldots,s-1)$, $A_s=1$. 
$($For the details of the choices of $A_i$, see $[30]$.$)$ Then the following identity holds$:$
\begin{equation*}
\begin{aligned}
&{}_{2s+3}F_{2s+2}\left(
\begin{array}{c}
a, \frac{a}{2}+1, c_0, b_1, c_1, \ldots, b_s, c_s\\
\frac{a}{2}, 1+a-c_0, 1+a-b_1, 1+a-c_1, \ldots, 1+a-b_s, 1+a-c_s
\end{array}
;1\right)\\
&=\frac{\Gamma(1+a-b_s)\Gamma(1+a-c_s)}{\Gamma(1+a)\Gamma(1+a-b_s-c_s)}
\sum_{l_1,\ldots,l_s\ge0}\frac{(b_1)_{l_1}(c_1)_{l_1}}{{l_1}!(1+a-c_0)_{l_1}}\\
&\times\prod_{i=2}^{s}\frac{(1+a-b_{i-1}-c_{i-1})_{l_i}(b_i)_{l_1+\cdots+l_i}(c_i)_{l_1+\cdots+l_i}}
{l_i!(1+a-b_{i-1})_{l_1+\cdots+l_i}(1+a-c_{i-1})_{l_1+\cdots+l_i}}.
\end{aligned}
\end{equation*}
\end{kr}
In \cite{kr}, C. Krattenthaler and T. Rivoal used their hypergeometric identities \cite[Proposition 1 (i) and (ii)]{kr} to give 
an alternative proof of Zudilin's identity \cite[Theorem 5]{zu}, 
which is an identity between certain very-well-poised hypergeometric series and certain multiple integrals 
related to the construction of $\mathbb{Q}$-linear forms in the values of the Riemann zeta function at 
positive integers. 
As C. Krattenthaler and T. Rivoal stated in \cite{kr}, their hypergeometric identities \cite[Proposition 1 (i) and (ii)]{kr} 
are non-terminating versions of a limiting case of 
a basic hypergeometric identity of G. E. Andrews \cite[Theorem 4]{an}.
\begin{remark}
For obvious examples of Theorem A, we can get the known identities
\begin{equation*}
2(1-2^{1-2s})\zeta(2s)=\zeta^{\star}(\underbrace{2,\ldots,2}_{s})
\end{equation*}
for any $s\in\mathbb{Z}_{\ge1}$ (Aoki--Ohno \cite[Theorem 1]{ao}, Zlobin \cite[Theorem 2]{zl2}, 
Vasil'ev \cite[Theorem]{va}; see also \cite{ako}, \cite{aow}, \cite{ikoo}, \cite{mu2}), 
and 
\begin{equation*}
2\zeta(2s+1)=\zeta^{\star}(1,\underbrace{2,\ldots,2}_{s})
\end{equation*}
for any $s\in\mathbb{Z}_{\ge1}$ (Zlobin \cite{zl2}, Vasil'ev \cite[Theorem]{va}; 
see also \cite[Examples (b)]{ow}, \cite{oz}).
\end{remark}
\section{Applications of the hypergeometric identities of C. Krattenthaler and T. Rivoal}
In this section, by using several ideas of other researchers, we derive several relations among MZSVs 
from the hypergeometric identities of C. Krattenthaler and T. Rivoal in Theorem A.
\par
(\textbf{A1}) Taking $a=2\alpha$, $b_1=1$, $b_i=\alpha$ ($i=2,\ldots,s+1$) and $c_j=\alpha$ ($j=1,\ldots,s+1$), 
where $s\in\mathbb{Z}_{\ge1}$ and $\alpha\in\mathbb{C}$ with $\mathrm{Re}\,\alpha>0$, 
in Theorem A (i), we get the identity
\begin{equation*}
\begin{aligned}
&\sum_{m=0}^{\infty}\frac{(-1)^m}{(m+\alpha)^{2s}}\\
=
&\frac{\Gamma(\alpha)^2}{2\Gamma(2\alpha)}
\sum_{0\le m_1\le \cdots{\le}m_s}\frac{(\alpha)_{m_1}^2}{{m_1}!(2\alpha)_{m_1}}
\frac{1}{(m_1+\alpha)(m_2+\alpha)^2 \cdots (m_s+\alpha)^2}.
\end{aligned}
\end{equation*}
(We note that this kind of identity can be found in \'{E}mery \cite[Proposition]{em}.) 
Differentiating both sides of this identity at $\alpha=1$, 
and using the equalities
\begin{equation*}
\begin{aligned}
&\frac{\mathrm{d}}{\mathrm{d}\alpha}\left(\frac{\Gamma(\alpha)^2}{2\Gamma(2\alpha)}\right){\Bigl|}_{\alpha=1}
=
-1,\\
&\frac{\mathrm{d}}{\mathrm{d}\alpha}(\alpha)_{m_1}{\Bigl|}_{\alpha=1}
=
{m_1}!\left(\sum_{l=0}^{m_1}\frac{1}{l+1}-\frac{1}{m_1+1}\right),\\
&\frac{\mathrm{d}}{\mathrm{d}\alpha}\frac{1}{(2\alpha)_{m_1}}{\Bigl|}_{\alpha=1}
=
\frac{-2}{(m_1+1)!}\left(\sum_{l=0}^{m_1}\frac{1}{l+1}-1\right),\\
\end{aligned}
\end{equation*}
we get
\begin{equation*}
\begin{aligned}
&-2s\sum_{m=0}^{\infty}\frac{(-1)^m}{(m+1)^{2s+1}}\\
=&-\zeta^{\star}(\underbrace{2,\ldots,2}_{s})+\zeta^{\star}(1,\underbrace{2,\ldots,2}_{s})
-\zeta^{\star}(3,\underbrace{2,\ldots,2}_{s-1})\\
&-\zeta^{\star}(1,\underbrace{2,\ldots,2}_{s})+\zeta^{\star}(\underbrace{2,\ldots,2}_{s})\\
&-\frac{1}{2}\zeta^{\star}(3,\underbrace{2,\ldots,2}_{s-1})
-\sum_{i=1}^{s-1}\zeta^{\star}(\underbrace{2,\ldots,2}_{i},3,\underbrace{2,\ldots,2}_{s-1-i})\\
=&-\frac{1}{2}\zeta^{\star}(3,\underbrace{2,\ldots,2}_{s-1})
-\sum_{i=1}^{s}\zeta^{\star}(\underbrace{2,\ldots,2}_{i-1},3,\underbrace{2,\ldots,2}_{s-i})
\end{aligned}
\end{equation*}
for any $s\in\mathbb{Z}_{\ge1}$. 
Since the identity
\begin{equation*}
\sum_{m=0}^{\infty}\frac{(-1)^m}{(m+1)^{2s+1}}
=
(1-2^{-2s})\zeta(2s+1)
\end{equation*}
holds for any $s\in\mathbb{Z}_{\ge1}$, we get the identity for $\zeta(2s+1)$ of K. Ihara, J. Kajikawa, Y. Ohno and J. Okuda 
in \cite[Theorem 2]{ikoo}:
\begin{equation}
\begin{aligned}
&4s(1-2^{-2s})\zeta(2s+1)\\
=
&\zeta^{\star}(3,\underbrace{2,\ldots,2}_{s-1})
+2\sum_{i=1}^{s}\zeta^{\star}(\underbrace{2,\ldots,2}_{i-1},3,\underbrace{2,\ldots,2}_{s-i})
\end{aligned}
\end{equation}
for any $s\in\mathbb{Z}_{\ge1}$. 
In \cite{ikoo}, K. Ihara, J. Kajikawa, Y. Ohno and J. Okuda proved the identity (1) 
by using the derivation relation for MZVs \cite[Corollary 6]{ikz}, 
and the relation among MZSVs of T. Aoki and Y. Ohno in \cite[Theorem 1]{ao}.
\begin{chu}[November 14, 2011]
To the best of my memory, the identity for $\zeta(2s+1)$ of K. Ihara, J. Kajikawa, Y. Ohno and J. Okuda 
in \cite[Theorem 2]{ikoo} was stated by Yasuo Ohno in his talk at 
the workshop ``Zeta Wakate Kenky\={u}sh\={u}kai'', Nagoya University, Japan, February 17--18, 2007: 
Yasuo Ohno talked on his collaboration with K. Ihara and J. Okuda.
\end{chu}
\par
(\textbf{A2}) Taking $a=2\alpha$, $c_0=1$ and $b_i=c_i=\alpha$ ($i=1,\ldots,s$), 
where $s\in\mathbb{Z}_{\ge2}$ and $\alpha\in\mathbb{C}$ with $\mathrm{Re}\,\alpha>0$, 
in Theorem A (ii), we get the identity 
\begin{equation*}
\begin{aligned}
&\sum_{m=0}^{\infty}\frac{1}{(m+\alpha)^{2s-1}}\\
=
&\frac{\Gamma(\alpha)^2}{2\Gamma(2\alpha)}
\sum_{0\le m_1\le \cdots{\le}m_s}\frac{(\alpha)_{m_1}^2}{{m_1}!(2\alpha)_{m_1}}
\frac{1}{(m_2+\alpha)^2 \cdots (m_s+\alpha)^2}.
\end{aligned}
\end{equation*}
(We note that this kind of identity can be found in \'{E}mery \cite[Proposition]{em}.) 
Differentiating both sides of this identity at $\alpha=1$, 
and by the same calculation as in (\textbf{A1}), 
we can get the identity
\begin{equation*}
\begin{aligned}
&(2s-1)\zeta(2s)\\
=&\zeta^{\star}(\underbrace{2,\ldots,2}_{s})
+\sum_{i=0}^{s-2}\zeta^{\star}(1,\underbrace{2,\ldots,2}_{i},3,\underbrace{2,\ldots,2}_{s-2-i})
\end{aligned}
\end{equation*}
for any $s\in\mathbb{Z}_{\ge1}$. This is an example of the cyclic sum formula for MZSVs 
which was proved by Y. Ohno and N. Wakabayashi in \cite[Theorem 1, Examples (a)]{ow}.
\par
(\textbf{A3}) Taking $a=2$, $c_0=\alpha$ and $b_i=c_i=1$ ($i=1,\ldots,s$), 
where $s\in\mathbb{Z}_{\ge2}$ and $\alpha\in\mathbb{C}$ with $\mathrm{Re}\,\alpha<2$, 
in Theorem A (ii), and multiplying both sides of the result by $(2-\alpha)^{-1}$, we get the identity
\begin{equation*}
\begin{aligned}
&\sum_{m=0}^{\infty}\frac{(\alpha)_m}{(2-\alpha)_{m+1}}\frac{1}{(m+1)^{2s-2}}\\
=
&\frac{1}{2}
\sum_{0\le m_1\le \cdots{\le}m_s}\frac{{m_1}!}{(2-\alpha)_{m_1+1}}
\frac{1}{(m_2+1)^2 \cdots (m_s+1)^2}.
\end{aligned}
\end{equation*}
By using an idea used by G. Kawashima in \cite{kaw} (i.e., making use of the product of finite multiple harmonic sums 
to derive relations among multiple zeta(-star) values (see also Vermaseren \cite[Section 5, Appendices A--F]{ve})), 
we can derive a relation among MZSVs from the above hypergeometric identity. 
Indeed, differentiating both sides of the above identity $r$ times at $\alpha=1$, 
and using the equalities
\begin{align}
\nonumber&\frac{1}{r!}\frac{\mathrm{d}^r}{\mathrm{d}\alpha^r}\left(\frac{(\alpha)_m}{(2-\alpha)_{m+1}}\right){\Bigl|}_{\alpha=1}\\
=
\nonumber&\frac{1}{m+1}\sum_{i=0}^{r}S_m(\underbrace{1,\ldots,1}_{r-i})
S^{\star}_m(\underbrace{1,\ldots,1}_{i}),\\
\nonumber&\frac{1}{r!}\frac{\mathrm{d}^r}{\mathrm{d}\alpha^r}\frac{1}{(2-\alpha)_{m+1}}{\Bigl|}_{\alpha=1}\\
=
&\frac{1}{(m+1)!}S^{\star}_m(\underbrace{1,\ldots,1}_{r})
\end{align}
for any $r\in\mathbb{Z}_{\ge0}$, where
\begin{equation*}
\begin{aligned}
&S_m(k_1,\ldots,k_n) := \sum_{0\le m_1< \cdots <m_n<m}\frac{1}{(m_1+1)^{k_1} \cdots (m_n+1)^{k_n}},\\
&S^{\star}_m(k_1,\ldots,k_n) := \sum_{0\le m_1\le \cdots \le m_n\le m}\frac{1}{(m_1+1)^{k_1} \cdots (m_n+1)^{k_n}},\\
&S_m(\underbrace{k_1,\ldots,k_n}_{0})=S^{\star}_m(\underbrace{k_1,\ldots,k_n}_{0})=1
\end{aligned}
\end{equation*}
for $k_1,\ldots,k_n\in\mathbb{Z}_{\ge1}$, $m\in\mathbb{Z}_{\ge0}$, we get the identity
\begin{equation}
\begin{aligned}
&\zeta^{\star}(\underbrace{1,\ldots,1}_{r+1},\underbrace{2,\ldots,2}_{s-1})\\
=
&2\sum_{m=0}^{\infty}\frac{1}{(m+1)^{2s-1}}
\left(\sum_{i=0}^{r}S_m(\underbrace{1,\ldots,1}_{r-i})
S^{\star}_m(\underbrace{1,\ldots,1}_{i})\right)
\end{aligned}
\end{equation}
for any $r\in\mathbb{Z}_{\ge0}$, $s\in\mathbb{Z}_{\ge2}$. 
We note that an equality like the equality (2) was already used by M. E. Hoffman in \cite[Proof of Corollary 4.2]{ho} 
to prove a relation among multiple zeta values: he used the derivative of the binomial coefficient 
(see also Kawashima \cite[Proofs of Propositions 4.7 and 5.2]{kaw}). 
By using the same method as used by G. Kawashima in \cite{kaw} (see also Vermaseren \cite{ve}), 
we see that the identity (3) is a relation among MZSVs. 
For example, by calculating the products of finite multiple harmonic sums on the right-hand side of the identity (3), 
the cases $r=0, 1, 2, 3$ in the identity (3) become as follows:
\begin{equation*}
\begin{aligned}
\zeta^{\star}(1,\underbrace{2,\ldots,2}_{s-1})
=&2\zeta(2s-1),\\
\zeta^{\star}(1,1,\underbrace{2,\ldots,2}_{s-1})
=&2^2\zeta^{\star}(1,2s-1)-2\zeta(2s),\\
\zeta^{\star}(1,1,1,\underbrace{2,\ldots,2}_{s-1})
=&2^3\zeta^{\star}(1,1,2s-1)-2^2\zeta^{\star}(2,2s-1)\\
&-2^2\zeta^{\star}(1,2s)+2\zeta(2s+1),\\
\zeta^{\star}(1,1,1,1,\underbrace{2,\ldots,2}_{s-1})
=&2^4\zeta^{\star}(1,1,1,2s-1)-2^3\zeta^{\star}(2,1,2s-1)\\
&-2^3\zeta^{\star}(1,2,2s-1)-2^3\zeta^{\star}(1,1,2s)\\
&+2^2\zeta^{\star}(3,2s-1)+2^2\zeta^{\star}(2,2s)\\
&+2^2\zeta^{\star}(1,2s+1)-2\zeta(2s+2)
\end{aligned}
\end{equation*}
for any $s\in\mathbb{Z}_{\ge2}$. 
The first identity above was already stated in Remark 1. 
The second identity above was proved by Y. Ohno and W. Zudilin in \cite[Lemma 5]{oz}. 
By the above examples of the identity (3), we observe that the right-hand side of the identity (3) probably coincides with 
the expression of the two-one formula for $\zeta^{\star}(\underbrace{1,\ldots,1}_{s_1},\underbrace{2,\ldots,2}_{s_2})$ 
($s_1, s_2\in\mathbb{Z}_{\ge1}$), i.e., the right-hand side of the identity (7a) in \cite{oz} for 
the index $(\underbrace{1,\ldots,1}_{s_1},\underbrace{2,\ldots,2}_{s_2})$ ($s_1, s_2\in\mathbb{Z}_{\ge1}$). 
The two-one formula for MZSVs was conjectured by Y. Ohno and W. Zudilin in \cite[pp.~327--328]{oz}. 
They proved two special cases of the two-one formula for MZSVs in \cite[Theorems 1 and 2]{oz}. 
In particular, they proved the two-one formula for $\zeta^{\star}(\underbrace{1,\ldots,1}_{s},2)$ ($s\in\mathbb{Z}_{\ge1}$) 
in \cite[Theorem 2]{oz}.
\par
(\textbf{A4}) Taking $a=2$, $b_1=\alpha$, $b_i=1$ ($i=2,\ldots,s+1$) and $c_j=1$ ($j=1,\ldots,s+1$), 
where $s\in\mathbb{Z}_{\ge1}$ and $\alpha\in\mathbb{C}$ with $\mathrm{Re}\,\alpha<3/2$, 
in Theorem A (i), we get the identity
\begin{equation*}
\begin{aligned}
&\sum_{m=0}^{\infty}\frac{(\alpha)_m}{(2-\alpha)_{m+1}}\frac{(-1)^m}{(m+1)^{2s-1}}\\
=
&\frac{1}{2}
\sum_{0\le m_1\le \cdots{\le}m_s}
\frac{1}{(m_1+2-\alpha)(m_1+1)(m_2+1)^2 \cdots (m_s+1)^2}.
\end{aligned}
\end{equation*}
We note that this kind of identity was already used by M. E. Hoffman in \cite[Section 4]{ho} to prove 
a relation among multiple zeta values: he used an identity of L. J. Mordell. 
By the same calculation as in (\textbf{A3}), we can get the identity
\begin{equation}
\begin{aligned}
&\zeta^{\star}(r+2,\underbrace{2,\ldots,2}_{s-1})\\
=
&2\sum_{m=0}^{\infty}\frac{(-1)^m}{(m+1)^{2s}}
\left(\sum_{i=0}^{r}S_m(\underbrace{1,\ldots,1}_{r-i})
S^{\star}_m(\underbrace{1,\ldots,1}_{i})\right)
\end{aligned}
\end{equation}
for any $r\in\mathbb{Z}_{\ge0}$, $s\in\mathbb{Z}_{\ge1}$. 
By the same calculation for the product of finite multiple harmonic sums as in (\textbf{A3}) (cf. Kawashima \cite{kaw}, Vermaseren \cite{ve}), 
the cases $r=0, 1, 2, 3$ in the identity (4) become as follows:
\begin{equation*}
\begin{aligned}
\zeta^{\star}(\underbrace{2,\ldots,2}_{s})
=&2\zeta_{-}^{\star}(2s),\\
\zeta^{\star}(3,\underbrace{2,\ldots,2}_{s-1})
=&2^2\zeta_{-}^{\star}(1,2s)-2\zeta_{-}^{\star}(2s+1),\\
\zeta^{\star}(4,\underbrace{2,\ldots,2}_{s-1})
=&2^3\zeta_{-}^{\star}(1,1,2s)-2^2\zeta_{-}^{\star}(2,2s)\\
&-2^2\zeta_{-}^{\star}(1,2s+1)+2\zeta_{-}^{\star}(2s+2),\\
\zeta^{\star}(5,\underbrace{2,\ldots,2}_{s-1})
=&2^4\zeta_{-}^{\star}(1,1,1,2s)-2^3\zeta_{-}^{\star}(2,1,2s)\\
&-2^3\zeta_{-}^{\star}(1,2,2s)-2^3\zeta_{-}^{\star}(1,1,2s+1)\\
&+2^2\zeta_{-}^{\star}(3,2s)+2^2\zeta_{-}^{\star}(2,2s+1)\\
&+2^2\zeta_{-}^{\star}(1,2s+2)-2\zeta_{-}^{\star}(2s+3)
\end{aligned}
\end{equation*}
for any $s\in\mathbb{Z}_{\ge1}$, where
\begin{equation*}
\zeta_{-}^{\star}(k_1,\ldots,k_n)
:=
\sum_{0< m_1\le \cdots{\le}m_n}
\frac{(-1)^{m_n-1}}{m_1^{k_1} \cdots m_n^{k_n}}
\end{equation*}
for $k_1,\ldots,k_n\in\mathbb{Z}_{\ge1}$: this alternating multiple series is a special case of 
Euler sums (see, e.g., \cite{bbv}, \cite{bbg}, \cite{bbb}, \cite{bowbra}, \cite{b}, \cite{bk}, \cite{eu}, \cite{fs}, \cite{ma}, 
\cite{n}, \cite{s}, \cite{ss}, \cite{zh}) 
or a special value of multiple polylogarithms (see, e.g., \cite{bbbl}, \cite{bb}, \cite{g}, \cite{g2}, 
\cite{ra}, \cite{u}, \cite{zh2}, \cite{zl3}). 
(For bases of the $\mathbb{Q}$-vector space generated by all Euler sums, 
see, e.g., Bl\"{u}mlein--Broadhurst--Vermaseren \cite{bbv}, 
Brown \cite[Subsection 5.3. Comments ii) and v)]{brown}, Deligne \cite[Section 7]{de}, Zlobin \cite[Section 8]{zl3}.) 
The first identity among the above examples of the identity (4) was already stated in Remark 1. 
By the above examples of the identity (4), we observe that the multiple zeta-star values 
$\zeta^{\star}(s_1+2,\underbrace{2,\ldots,2}_{s_2})$ ($s_1, s_2\in\mathbb{Z}_{\ge0}$) and the Euler sums 
$\zeta_{-}^{\star}(k_1,\ldots,k_n)$ probably satisfy a relation like the two-one formula for 
$\zeta^{\star}(\underbrace{1,\ldots,1}_{s_1},\underbrace{2,\ldots,2}_{s_2})$ ($s_1, s_2\in\mathbb{Z}_{\ge1}$). 
We note that relations between multiple zeta values and Euler sums similar to the above examples of the identity (4) were already studied by 
D. J. Broadhurst in \cite[p.~14, Eq. (32)]{b}, and J. Bl\"{u}mlein, D. J. Broadhurst and J. A. M. Vermaseren in \cite[Sections 2 and 4]{bbv}: 
they called the relations the doubling relations. (See also Borwein--Bradley--Broadhurst--Lison\v{e}k \cite[Section 9]{bbbl}, 
Vermaseren \cite[Section 5]{ve}, Goncharov \cite[Proposition 2.3]{g}, \cite[Proposition 2.13]{g2}.) 
We also note the results for multiple $L$-values of M. Nishi which were stated in \cite[Proposition 4.2]{ak}.
\begin{remark}
In order to derive the relations among MZSVs in (\textbf{A1})--(\textbf{A4}) from 
the hypergeometric identities of C. Krattenthaler and T. Rivoal in \cite[Proposition 1 (i) and (ii)]{kr}, 
we expressed hypergeometric series as $\mathbb{Q}$-linear combinations of MZSVs. 
This kind of expression was already studied in, e.g., \cite{A1}, \cite{A2}, \cite{A3}, \cite{A4}, \cite{A5}, \cite{A6} 
to evaluate Feynman diagrams; in, e.g., \cite{ako}, \cite{ao}, \cite{aow}, \cite{A7}, \cite{A8} to obtain relations among 
multiple zeta(-star) values; and in, e.g., \cite{cfr}, \cite{cfr2}, \cite{fi}, \cite{zl} 
to construct $\mathbb{Q}$-linear forms in multiple zeta(-star) values.
\end{remark}
Though the contents of this note depend on many prior works and several ideas of other researchers, 
the hypergeometric identities of C. Krattenthaler and T. Rivoal in 
\cite[Proposition 1 (i) and (ii)]{kr} seem to be useful for the study of relations among multiple zeta(-star) values 
(see also \cite[Remark 2.7]{ig}).
\section*{Addendum (July 12, 2011)}
We give proofs for the observations in (\textbf{A3}) and (\textbf{A4}). 
By direct calculation, we can easily see that the identities
\begin{equation*}
\begin{aligned}
\frac{(\alpha)_m}{(2-\alpha)_{m+1}}
=&\frac{1}{2-\alpha+m}\prod_{i=0}^{m-1}\left(1+\frac{2(\alpha-1)}{2-\alpha+i}\right)\\
=&\frac{1}{2-\alpha+m}\sum_{i=0}^{m}2^i\sum_{0\le m_1<\cdots<m_i<m}\prod_{j=1}^{i}\left(\frac{\alpha-1}{2-\alpha+m_j}\right)
\end{aligned}
\end{equation*}
hold for any $m\in\mathbb{Z}_{\ge0}$. 
Using the second identity above, we can prove the equality
\begin{equation}
\begin{aligned}
&\frac{1}{r!}\frac{\mathrm{d}^r}{\mathrm{d}\alpha^r}\left(\frac{(\alpha)_m}{(2-\alpha)_{m+1}}\right){\Bigl|}_{\alpha=1}\\
=&\sum_{i=0}^{m}2^i
\sum_{\begin{subarray}{c}
k_1+\cdots+k_{i+1}=r+1\\
k_j\in\mathbb{Z}_{\ge1}
\end{subarray}}
S_m(k_1,\ldots,k_i)\frac{1}{(m+1)^{k_{i+1}}}
\end{aligned}
\end{equation}
for any $r, m\in\mathbb{Z}_{\ge0}$. 
The equality (5) gives us the identity
\begin{equation*}
\begin{aligned}
&2\sum_{m=0}^{\infty}\frac{1}{(m+1)^{2s-2}}
\left(\frac{1}{r!}\frac{\mathrm{d}^r}{\mathrm{d}\alpha^r}\left(\frac{(\alpha)_m}{(2-\alpha)_{m+1}}\right){\Bigl|}_{\alpha=1}\right)\\
=&\sum_{i=0}^{r}2^{i+1}
\sum_{\begin{subarray}{c}
k_1+\cdots+k_{i+1}=r+1\\
k_j\in\mathbb{Z}_{\ge1}
\end{subarray}}
\zeta(k_1,\ldots,k_i,k_{i+1}+2s-2)
\end{aligned}
\end{equation*}
for any $r\in\mathbb{Z}_{\ge0}$, $s\in\mathbb{Z}_{\ge2}$, where
\begin{equation*}
\zeta(k_1,\ldots,k_n)
:=
\sum_{0< m_1< \cdots < m_n}
\frac{1}{m_1^{k_1} \cdots m_n^{k_n}}
\end{equation*}
($k_1,\ldots,k_n\in\mathbb{Z}_{\ge1}$, $k_n\ge2$) is the multiple zeta value. 
By using this identity and the argument used by Y. Ohno and W. Zudilin in 
\cite[p.~328, Proof of the equality of the right-hand sides in (7a) and (7b)]{oz}, 
the identity (3) can be rewritten as
\begin{equation*}
\begin{aligned}
&\zeta^{\star}(\underbrace{1,\ldots,1}_{r+1},\underbrace{2,\ldots,2}_{s-1})\\
=
&\sum_{i=0}^{r}(-1)^{r-i}2^{i+1}
\sum_{\begin{subarray}{c}
k_1+\cdots+k_{i+1}=r+1\\
k_j\in\mathbb{Z}_{\ge1}
\end{subarray}}
\zeta^{\star}(k_1,\ldots,k_i,k_{i+1}+2s-2)
\end{aligned}
\end{equation*}
for any $r\in\mathbb{Z}_{\ge0}$, $s\in\mathbb{Z}_{\ge2}$. 
This is the two-one formula for 
$\zeta^{\star}(\underbrace{1,\ldots,1}_{s_1},\underbrace{2,\ldots,2}_{s_2})$ ($s_1, s_2\in\mathbb{Z}_{\ge1}$).
\par
Similarly, by using the equality (5) and the argument used by Y. Ohno and W. Zudilin in 
\cite[p.~328, Proof of the equality of the right-hand sides in (7a) and (7b)]{oz}, 
the identity (4) can be also rewritten as
\begin{equation*}
\begin{aligned}
&\zeta^{\star}(r+2, \underbrace{2,\ldots,2}_{s-1})\\
=
&\sum_{i=0}^{r}(-1)^{r-i}2^{i+1}
\sum_{\begin{subarray}{c}
k_1+\cdots+k_{i+1}=r+1\\
k_j\in\mathbb{Z}_{\ge1}
\end{subarray}}
\zeta_{-}^{\star}(k_1,\ldots,k_i,k_{i+1}+2s-1)
\end{aligned}
\end{equation*}
for any $r\in\mathbb{Z}_{\ge0}$, $s\in\mathbb{Z}_{\ge1}$.
\begin{flushleft}
Graduate School of Mathematics\\
Nagoya University\\
Chikusa-ku, Nagoya 464-8602, Japan\\
\textit{E-mail address}: m05003x@math.nagoya-u.ac.jp\\
\end{flushleft}
\end{document}